\documentclass[11pt]{article}
\usepackage{amscd}
\usepackage{amsfonts}
\usepackage{amsmath}
\usepackage{amssymb}
\usepackage{amsthm}
\usepackage{bbm}
\usepackage{CJK}
\usepackage{fancyhdr}
\usepackage{graphicx}
\usepackage{indentfirst}
\usepackage{latexsym}
\usepackage{mathrsfs}
\usepackage{xypic}

\allowdisplaybreaks[4]
\usepackage[top=1in,bottom=1in,left=1.25in,right=1.25in]{geometry}
\def\<{\langle}
\def\>{\rangle}

\def\o{\otimes}

\date{}
\begin{document}
\renewcommand{\baselinestretch}{1.2}
\renewcommand{\arraystretch}{1.0}
\title{\bf    Drinfel'd double for monoidal Hom-Hopf algebras}
\author{{\bf Yan Ning$^{1,2}$, Daowei Lu$^{1,2}$\footnote {Corresponding author: Email: ludaowei620@126.com }, Xiaohui Zhang$^1$}\\
{\small $^1$School of Mathematical Sciences, Qufu Normal University}\\
{\small Qufu, Shandong 273165, P. R. China}\\
{\small $^2$Department of Mathematics, Jining University}\\
{\small Qufu, Shandong 273155, P. R. China}
}
\maketitle
\begin{center}
\begin{minipage}{12.cm}

\noindent{\bf Abstract.} In this paper we mainly  construct bicrossproduct for finite-dimensional monoidal Hom-Hopf algebra $(H,\alpha)$, generalizing the Majid's bicrossproduct. Naturally the Hom-type bicrossproduct  leads to Drinfel'd double $(H^{op}\bowtie H^{\ast},\alpha\otimes(\alpha^{-1})^{\ast})$  with a quasitriangular structure $R$  satisfying
   the quantum Hom-Yang-Baxter equations.
 \\

\\

 \end{minipage}
 \end{center}
 \normalsize\vskip1cm

\section*{INTRODUCTION}

Algebraic deformation has been well developed recently into a new broad class of non-associative algebras, and its theory has been applied in modules of quantum phenomena, as well as in analysis
of complex systems.
In \cite{HLS} by the help of general $\sigma$-derivations, Hartwig, Larsson and Silvestrov developed a general quasi-deformation scheme for
Lie algebras of vector fields. Then the deformation has been investigated by Makhlouf and
Silvestrov in \cite{MS1}. In their construction of Hom-Lie algebra, the Jacobi identity is replaced by the so-called Hom-Jacobi identity, namely
$$[\alpha(x),[y,z]]+[\alpha(y),[z,x]]+[\alpha(z),[x,y]]=0,$$
where $\alpha$ is an endomorphism of the Lie algebra. The main initial goal of this investigation was to create a
unified general approach to examples of $q$-deformations of Witt and Virasoro algebras. Particularly it
was observed that in these examples some $q$-deformations of ordinary Lie algebra Jacobi
identities hold. Motivated by these new interesting examples, quasi-Hom-Lie algebra and Hom-Lie algebra were introduced in \cite{HLS}.

The idea of construction by adapting associativity-like conditions via endomorphisms was applied to other
algebraic structures. For example,  Hom-algebras and Hom-coalgebra was introduced for the first time in \cite{MS1} and \cite{MS2}, respectively. Here the Hom-associativity and Hom-coassociativity were obtained by twisting the original associativity and coassociativity by endomorphisms. Hom-module and Hom-comodule were naturally developed in the same way. Then Hom-bialgebras were defined in a rational way.  These objects are slightly different from what we will study (called monoidal Hom-Hopf algebra) in this paper, as introduced in \cite{CG}.

 In \cite{LS} Liu and Shen constructed  Radford's biproduct on monoidal Hom-Hopf algebras. In \cite{Y} D. Yau introduced the quasitriangular Hom-bialgebras(not monoidal Hom-bialgebras), as a generalization of the ordinary quasitriangular bialgebras and the quantum Hom-Yang-Baxter equation(QHYBE) of the form
\begin{align*}
&R^{12}(R^{13}R^{23})=(R^{13}R^{23})R^{12}\\
&(R^{12}R^{13})R^{23}=R^{23}(R^{13}R^{12}).
\end{align*}
The construction of Majid's bicrossproduct Hopf algebras was motivated by the search for examples of self-dual algebraic structures, which means in the first place to find a category with a dualising endofunctor such that some kind of Pontryagin duality theorem holds. The bicrossproduct Hopf algebras provide numerous examples of non-commutative and non-cocomutative Hopf algebras, and moreover they turn out be closely related to the Drinfel'd double.

  Motivated by these ideas, in this paper as a continuation of \cite{LS}, we will construct the Majid's bicrossproduct for monoidal Hom-Hopf algebras and then in the framework of monoidal Hom-Hopf algebras,
we will consider the Drinfel'd double which could be deduced from the Majid's bicrossproduct Hopf algebra and the its quasitriangular structure.

This paper is organized as follows.

In Section 1, we will recall the definitions and results of monoidal Hom-Hopf algebras, such as a monoidal Hom-algeba, a monoidal Hom-coalgebra, a  monoidal Hom-module, a monoidal Hom-comodule and the Hom-smash products.

In Section 2, we will introduce the notion of bicrossproduct $(B\#H,\beta\#\alpha)$, and give the necessary and sufficient conditions for $(B\#H,\beta\#\alpha)$ to form a monoidal Hom-Hopf algebra (see Theorem 2.2 and Theorem 2.3), generalizing the Majid's bicrossproduct defined in \cite{M1}. An example will be given at the end of this section.

In Section 3, we will construct a class of bicrossproduct monoidal Hom-Hopf algebras (see Theorem 3.3).

In Section 4, we will construct the generalized Drinfel'd double
 from the Majid's bicrossproduct for monoidal Hopf algebras (see Theorem 4.2) and obtain its quasitriangular structure
 (see Proposition 4.10), providing a solution of quantum Hom-Yang-Baxter equation.

Throughout this article, all the vector spaces, tensor product and homomorphisms are over a fixed field $k$. We use the Sweedler's notation for the terminologies on coalgebras. For a coalgebra $C$, we write comultiplication $\Delta(c)=\sum c_{1}\otimes c_{2}$ for any $c\in C$.

\section{PRELIMINARIES}

Let $\mathcal{M}_{k}=(\mathcal{M}_{k},\otimes,k,a,l,r)$ be the category of $k$-modules. Now from this category, we could construct a new monoidal category $\mathcal{H}(\mathcal{M}_{k})$. The objects of $\mathcal{H}(\mathcal{M}_{k})$ are pairs $(M,\mu)$, where $M\in\mathcal{M}_{k}$ and $\mu\in Aut_{k}(M)$. Any morphism $f:(M,\mu)\rightarrow (N,\nu)$ in $\mathcal{H}(\mathcal{M}_{k})$ is a $k$-linear map from $M$ to $N$ such that $\nu\circ f=f\circ \mu$. For any objects $(M,\mu)$ and $(N,\nu)$ in $\mathcal{H}(\mathcal{M}_{k})$, the monoidal structure is given by
$$(M,\mu)\otimes(N,\nu)=(M\otimes N,\mu\otimes\nu),$$
and the unit is $(k,id_{k})$.

Generally speaking, all Hom-structure are objects in the monoidal category $\tilde{\mathcal{H}}(\mathcal{M}_{k})=(\mathcal{H}(\mathcal{M}_{k}),\otimes,(k,id_{k}),\tilde{a},\tilde{l},\tilde{r})$ as introduced in \cite{CG}, where the associativity constraint $\tilde{a}$ is given by the formula
$$\tilde{a}_{M,N,L}=a_{M,N,L}\circ((\mu\otimes id)\otimes\lambda^{-1})=(\mu\otimes(id\otimes\lambda^{-1}))\circ a_{M,N,L}$$
for any objects $(M,\mu),\ (N,\nu),\ (L,\lambda)$ in $\mathcal{H}(\mathcal{M}_{k}).$ And the unit constraints $\tilde{l}$ and $\tilde{r}$ are defined by
$$\tilde{l}_{M}=\mu\circ l_{M}=l_{M}\circ(id\otimes\mu),\ \tilde{r}_{M}=\mu\circ r_{M}=r_{M}\circ(\mu\otimes id).$$
The category $\tilde{\mathcal{H}}(\mathcal{M}_{k})$ is called the Hom-category associated to the monoidal category $\mathcal{M}_{k}$.

In what follows, we will recall the definitions in \cite{CG} on the monoidal Hom-associative algebras, monoidal Hom-coassociative coalgebras, monoidal Hom-modules and monoidal Hom-comodules.

 A unital monoidal Hom-associative algebra is an object $(A,\alpha)$ in the category $\tilde{\mathcal{H}}(\mathcal{M}_{k})$ together with an element $1_{A}\in A$ and a linear map $m:A\otimes A\rightarrow A,\ a\otimes b\mapsto ab$ such that
\begin{align*}
&\alpha(a)(bc)=(ab)\alpha(c),\ a1_{A}=\alpha(a)=1_{A}a,\\
&\alpha(ab)=\alpha(a)\alpha(b),\ \alpha(1_{A})=1_{A},
\end{align*}
for all $a,b,c\in A.$

In the setting of Hopf algebras, $m$ is called the Hom-multiplication, $\alpha$ is the twisting automorphism, and $1_{A}$ is the unit.
Let $(A,\alpha)$ and $(A',\alpha')$ be two monoidal Hom-algebras. A Hom-algebra map $f:(A,\alpha)\rightarrow(A',\alpha')$ is a linear map such that $f\circ\alpha=\alpha'\circ f$, $f(ab)=f(a)f(b)$ and $f(1_{A})=1_{A}$.  Obviously $(A^{op},\alpha)$ is also a Hom-algebra.

 A counital monoidal Hom-coassociative coalgebra is an object $(C,\gamma)$ in the category $\tilde{\mathcal{H}}(\mathcal{M}_{k})$ together with linear maps $\Delta:C\rightarrow C\otimes C,\ c\mapsto c_{1}\otimes c_{2}$ and $\varepsilon:C\rightarrow k$ such that
\begin{align*}
&\sum\gamma^{-1}(c_{1})\otimes\Delta(c_{2})=\sum\Delta(c_{1})\otimes\gamma^{-1}(c_{2}),\\
 &\sum c_{1}\varepsilon(c_{2})=\sum\varepsilon(c_{1})c_{2}=\lambda^{-1}(c),\\
&\Delta(\gamma(c))=\sum\gamma(c_{1})\otimes\gamma(c_{2}),\ \varepsilon\gamma(c)=\varepsilon(c),
\end{align*}
for all $c\in C.$

Let $(C,\gamma)$ and $(C',\gamma')$ be two monoidal Hom-coalgebras. A Hom-coalgebra map $f:(C,\gamma)\rightarrow(C',\gamma')$ is a linear map such that $f\circ\gamma=\gamma'\circ f,$ $\Delta\circ f=(f\otimes f)\circ\Delta$ and $\varepsilon\circ f=\varepsilon.$

A monoidal Hom-bialgebra $H=(H,\alpha,m,1_{H},\Delta,\varepsilon)$ is a bialgebra in the category $\tilde{\mathcal{H}}(\mathcal{M}_{k})$ if $(H,\alpha,m,1_{H})$ is a monoidal Hom-algebra and $(H,\alpha,\Delta,\varepsilon)$ is a monoidal Hom-coalgebra such that $\Delta$ and $\varepsilon$ are Hom-algebra maps, that is, for any $g,h\in H,$
\begin{align*}
&\Delta(gh)=\Delta(g)\Delta(h),\ \Delta(1_{H})=1_{H}\otimes 1_{H},\\
&\varepsilon(gh)=\varepsilon(g)\varepsilon(h),\ \varepsilon(1_{H})=1.
\end{align*}

A monoidal Hom-bialgebra $(H,\alpha)$ is called a monoidal Hom-Hopf algebra if there exists a linear map $S:H\rightarrow H$(the antipode) such that
$$S\circ\alpha=\alpha\circ S,\ \sum S(h_{1})h_{2}=\varepsilon(h)1_{H}=\sum h_{1}S(h_{2}).$$

Just as in the case of Hopf algebras, the antipode of monoidal Hom-Hopf algebras is a morphism of Hom-anti-algebras and Hom-anti-coalgebras.

When $(H,\alpha)$ is a finite-dimensional monoidal Hom-Hopf algebra, then $(H^{\ast},(\alpha^{-1})^{\ast})$ is also a monoidal Hom-Hopf algebra, where $(\alpha^{-1})^{\ast}:H^{\ast}\rightarrow H^{\ast},\ f\mapsto f\circ\alpha^{-1}$, and the antipode is $S^{\ast}_{H}$ defined similarly as $(\alpha^{-1})^{\ast}$.

Let $(A,\alpha)$ be a monoidal Hom-algebra. A left $(A,\alpha)$-Hom-module is an object $(M,\mu)$ in $\tilde{\mathcal{H}}(\mathcal{M}_{k})$ together with a linear map $\varphi:A\otimes M\rightarrow M,\ a\otimes m\mapsto am$ such that
$$
\alpha(a)(bm)=(ab)\mu(m),\ 1_{A}m=\mu(m),\ \mu(am)=\alpha(a)\mu(m),
$$
for all $a,b\in A$ and $m\in M$.

Similarly we can define the right $(A,\alpha)$-Hom-modules. Let $(M,\mu)$ and $(N,\nu)$ be two left $(A,\alpha)$-Hom-modules, then a linear map $f:M\rightarrow N$ is a called left $A$-module map if $f(am)=af(m)$ for any $a\in A$, $m\in M$ and $f\circ\mu=\nu\circ f$.

Let $(C,\gamma)$ be a monoidal Hom-coalgebra. A right $(C,\gamma)$-Hom-comodule is an object $(M,\mu)$ in $\tilde{\mathcal{H}}(\mathcal{M}_{k})$ together with a linear map $\rho_{M}:M\rightarrow M\otimes C,\ m\mapsto m_{(0)}\otimes m_{(1)}$ such that
\begin{align*}
&\sum\mu^{-1}(m_{(0)})\otimes\Delta(m_{(1)})=\sum\rho_{M}(m_{(0)})\otimes\gamma^{-1}(m_{(1)}),\\ &\sum\varepsilon(m_{(1)})m_{(0)}=\mu^{-1}(m),\\
&\rho_{M}(\mu(m))=\sum\mu(m_{(0)})\otimes\gamma(m_{(1)}),
\end{align*}
for all $m\in M.$

Let $(M,\mu)$ and $(N,\nu)$ be two right $(C,\gamma)$-Hom-comodules, then a linear map $g:M\rightarrow N$ is a called right $C$-comodule map if $g\circ \mu=\nu\circ g$ and $\rho_{N}(g(m))=(g\otimes id)\rho_{M}(m)$ for any $m\in M.$

Let $(H,\alpha)$ be a monoidal Hom-bialgebra, and $(B,\beta)$ be a monoidal Hom-algebra.
 We say $(B,\beta)$ is a left $(H,\alpha)$-Hom-module algebra  if $(B,\beta)$ is a left $(H,\alpha)$-Hom-module with the action $\cdot$ and satisfies
$$h\cdot(ab)=\sum(h_{1}\cdot a)(h_{2}\cdot b),\ h\cdot 1_{B}=\varepsilon(h)1_{B},$$
for any $a,b\in B$ and $h\in H$.

Dually,  let $(H,\alpha)$ be a monoidal Hom-bialgebra, and $(C,\gamma)$ be a monoidal Hom-coalgebra. We say $(C,\gamma)$ is a right $(H,\alpha)$-Hom-comodule coalgebra, if $(C,\gamma)$ is a right $(H,\alpha)$-Hom-comodule and with the coaction $\rho(c)=\sum c_{(0)}\otimes c_{(1)}$ and
\begin{enumerate}
\item
[(1)]$\sum c_{(0)1}\otimes c_{(0)2}\otimes c_{(1)}=\sum c_{1(0)}\otimes c_{2(0)}\otimes c_{1(1)}c_{2(1)},$
\item
[(2)]$\sum \varepsilon(c_{(0)})c_{1}=\varepsilon(c)1_{H}.$
\end{enumerate}

Let $(H,\alpha)$ be a monoidal Hom-bialgebra, and $(B,\beta)$ be a monoidal Hom-coalgebra. We say $(B,\beta)$ is a left $(H,\alpha)$-Hom-module coalgebra if $(B,\beta)$ is a left $(H,\alpha)$-Hom-module with the action $\cdot$ and satisfies
$$\Delta(h\cdot a)=\sum(h_{1}\cdot a_{1})(h_{2}\cdot a_{2}),\ \varepsilon_{B}(h\cdot a)=\varepsilon_{H}(h)\varepsilon_{B}(a),$$
for any $a\in B$ and $h\in H$.

Let $(B,\beta)$ be a left $(H,\alpha)$-Hom-module algebra. The Hom-smash product $(B\#H,\beta\#\alpha)$ of $(B,\beta)$ and $(H,\alpha)$ is defined as follows:
\begin{enumerate}
\item
[(1)]$B\#H=B\otimes H$ as a vector space,
\item
[(2)]For any $a,b\in B$ and $h,k\in H$
$$(a\#h)(b\#k)=\sum a(h_{1}\cdot\beta^{-1}(b))\#\alpha(h_{2})k.$$
\end{enumerate}
Then $(B\#H,\beta\#\alpha)$ is a monoidal Hom-associative algebra with the unit $1_{B}\#1_{H}$ (see \cite{LS}).

\section{The bicrossproduct construction}
\def\theequation{2.\arabic{equation}}
\setcounter{equation} {0}

In this section, we will construct the bicrossproduct of two monoidal Hom-Hopf algebras.
First, we give a slightly different version from the one in \cite{LS} for our purpose containing
 a simple sketch.

\noindent{\bf Proposition 2.1.}
Let $(H,\alpha)$ be a  monoidal Hom-coalgebra and $(B,\beta)$ be a monoidal Hom-bialgebra. Assume that $(H,\alpha)$ is a right  monoidal $(B,\beta)$-Hom comodule coalgebra, then $(B\#H,\beta\#\alpha)$ is a monoidal Hom-coalgebra with the comultiplication and counit defined as follows:
\begin{align}
&\Delta(a\#h)=\sum a_{1}\#\alpha(h_{1(0)})\otimes\beta^{-1}(a_{2})h_{1(1)}\# h_{2},\\
&\varepsilon(a\#h)=\varepsilon_{B}(a)\varepsilon_{H}(h),
\end{align}
where $B\#H=B\otimes H$ as a vector space, and we use the same notion as in the Hom-smash product.

\noindent{\bf Proof}~~Firstly, for any $a\#h\in B\#H$,
\begin{eqnarray*}
&&\sum\Delta((a\#h)_{1})\otimes(\beta^{-1}\#\alpha^{-1})((a\#h)_{2})\\
&=&\sum a_{11}\#\alpha(\alpha(h_{1(0)})_{1(0)})\otimes\beta^{-1}(a_{12})\alpha(h_{1(0)})_{1(1)}\#\alpha(h_{1(0)})_{2}\\
&&\otimes\beta^{-2}(a_{2})\beta^{-1}(h_{1(1)})\#\alpha^{-1}(h_{2})\\
&=&\sum\beta^{-1}(a_{1})\#\alpha(\alpha(\alpha^{-1}(h_{1})_{(0)})_{(0)})\otimes\beta^{-1}(a_{21})\alpha(\alpha^{-1}(h_{1})_{(0)})_{(1)}\\
&&\#\alpha(h_{21(0)})
\otimes\beta^{-1}(a_{22})\beta^{-1}(\alpha^{-1}(h_{1})_{(1)}h_{21(1)})\#h_{22} \ ( by \ coassociativity)\\
&=&\sum\beta^{-1}(a_{1})\#\alpha(\alpha(\alpha^{-1}(h_{1(0)}))_{(0)})\otimes\beta^{-1}(a_{21})\alpha(\alpha^{-1}(h_{1(0)}))_{(1)}\\
&&\#\alpha(h_{21(0)})\otimes\beta^{-1}(a_{22})\beta^{-1}(\beta^{-1}(h_{1(1)})h_{21(1)})\#h_{22}\\
&=&\sum\beta^{-1}(a_{1})\#\alpha(h_{1(0)(0)})\otimes\beta^{-1}(a_{21})h_{1(0)(1)})\\
&&\#\alpha(h_{21(0)})\otimes\beta^{-1}(a_{22})\beta^{-1}(\beta^{-1}(h_{1(1)})h_{21(1)})\#h_{22}\\
&=&\sum\beta^{-1}(a_{1})\#h_{1(0)}\otimes\beta^{-1}(a_{21})h_{1(1)1}\\
&&\#\alpha(h_{21(0)})\otimes\beta^{-1}(a_{22})\beta^{-1}(h_{1(1)2}h_{21(1)})\#h_{22}, \ (by\ comodule\ condition)\\
&=&\sum\beta^{-1}(a_{1})\# h_{1(0)}\otimes\beta^{-1}(a_{21})h_{1(1)1}\#\alpha(h_{21(0)})\\
&&\otimes(\beta^{-2}(a_{22})\beta^{-1}(h_{1(1)2}))h_{21(1)}\#h_{22}\\
&=&\sum(\beta^{-1}\#\alpha^{-1})(a_{1}\#\alpha(h_{1(0)}))\otimes\Delta(\beta^{-1}(a_{2})h_{1(1)}\#h_{2})\\
&=&\sum(\beta^{-1}\#\alpha^{-1})((a\#h)_{1})\otimes\Delta((a\#h)_{2})\\
\end{eqnarray*}
Hence we have the Hom-coassociativity.
It is straightforward to check that
\begin{eqnarray*}
&&\sum(\varepsilon_{B}\otimes\varepsilon_{H}\otimes id\otimes id)(a_{1}\#\alpha(h_{1(0)})\otimes\beta^{-1}(a_{2})h_{1(1)}\otimes h_{2})\\
&&\quad \quad \quad =
\sum(id\otimes id\otimes \varepsilon_{B}\otimes\varepsilon_{H})(a_{1}\#\alpha(h_{1(0)})\otimes\beta^{-1}(a_{2})h_{1(1)}\otimes h_{2}).
\end{eqnarray*}
This shows  that $(B\#H, \beta\#\alpha)$ is a monoidal Hom-coalgebra. The completes the proof. $\hfill \Box$
\\

Here we will call $(B\#H,\beta\#\alpha)$ a right Hom-smash coproduct as in \cite{LS}.

Before the next proposition, we have the following identities by the coassociativity:
\begin{align}
&\sum h_{11}\otimes h_{12}\otimes h_{211}\otimes h_{212}\otimes h_{22}=\sum\alpha^{-1}(h_{1})\otimes\alpha^{2}(h_{2111})\otimes \alpha(h_{2112})\otimes h_{212}\otimes h_{22},\\
&\sum h_{1}\otimes h_{211}\otimes h_{212}\otimes h_{22}=\sum\alpha(h_{11})\otimes \alpha^{-1}(h_{12})\otimes \alpha^{-1}(h_{21})\otimes h_{22}.
\end{align}

\noindent{\bf Theorem 2.2.}
Let $(B,\beta)$ and $(H,\alpha)$ be two monoidal Hom-Hopf algeberas. Let $(B\#H, \beta\#\alpha)$ be a Hom-smash product
 with a Hom-smash coproduct defined as above. Then
$(B\#H, \beta\#\alpha)$ is a Hom-bialgebra if and only if
\begin{enumerate}
\item
[(1)]$\Delta(h\cdot b)=\sum\alpha(h_{1(0)})\cdot b_{1}\otimes\beta(h_{1(1)})(\alpha^{-1}(h_{2})\cdot\beta^{-1}(b_{2}))$,
\item
[(2)]$\varepsilon_{B}(h\cdot b)=\varepsilon_{H}(h)\varepsilon_{B}(b),$
\item
[(3)]$\rho(1_{H})=1_{H}\otimes 1_{B}$,
\item
[(4)]$\sum h_{2(0)}\otimes(h_{1}\cdot b)\beta^{2}(h_{2(1)})=\sum h_{1(0)}\otimes\beta^{2}(h_{1(1)})(h_{2}\cdot b),$
\item
[(5)]$\sum(hk)_{(0)}\otimes(hk)_{(1)}=\sum\alpha(h_{1(0)})k_{(0)}\otimes\beta(h_{1(1)})(\alpha^{-1}(h_{2})\cdot\beta^{-1}(k_{(1)})).$
\end{enumerate}

\noindent{\bf Proof}~~If $(B\#H, \beta\#\alpha)$ is a monoidal Hom-bialgebra, by
$\Delta(1\#1)=1\#1\otimes 1\#1$, we have
$$
\sum1\#\alpha(1_{0})\otimes\beta(1_{(1)})\#1=1\#1\otimes 1\#1.
$$

Then
$$\sum1\#\alpha(1)_{0}\otimes\alpha(1)_{(1)}\#1=1\#1\otimes 1\#1,$$

hence
$$\sum1_{(0)}\otimes1_{(1)}=1_{H}\otimes 1_{B}.$$

And by $\varepsilon((a\#h)(b\#k))=\varepsilon(a\#h)\varepsilon(b\#k)$ we have
$$\sum\varepsilon(a(h_{1}\cdot\beta^{-1}(b))\#\alpha(h_{2})k)=\varepsilon_{B}(a)\varepsilon_{B}(b)\varepsilon_{H}(h)\varepsilon_{H}(k).$$

Then
$$\sum\varepsilon_{B}(a)\varepsilon_{B}(h_{1}\cdot\beta^{-1}(b))\varepsilon_{H}(h_{2}k)=\varepsilon_{B}(a)\varepsilon_{B}(b)\varepsilon_{H}(h)\varepsilon_{H}(k).$$

Hence
$$\varepsilon_{B}(h\cdot b)=\varepsilon_{H}(h)\varepsilon_{B}(b).$$

By a calculation as follows
\begin{align*}
&\Delta((a\#h)(b\#k))\\
=&\sum\Delta(a(h_{1}\cdot\beta^{-1}(b))\#\alpha(h_{2})k)\\
=&\sum a_{1}(h_{1}\cdot\beta^{-1}(b))_{1}\#\alpha((\alpha(h_{21})k_{1})_{(0)})\\
&\otimes\beta^{-1}(a_{2}(h_{1}\cdot\beta^{-1}(b))_{2})(\alpha(h_{21})k_{1})_{(1)}\#\alpha(h_{22})k_{2},
\end{align*}
and
\begin{align*}
&\Delta(a\#h)\Delta(b\#k)\\
=&\sum[a_{1}\#\alpha(h_{1(0)})\otimes\beta^{-1}(a_{2})h_{1(1)}\# h_{2}][b_{1}\#\alpha(k_{1(0)})\otimes\beta^{-1}(b_{2})k_{1(1)}\# k_{2}]\\
=&\sum a_{1}(\alpha(h_{1(0)})_{1}\cdot\beta^{-1}(b_{1}))\#\alpha(\alpha(h_{1(0)})_{2})\alpha(k_{1(0)})\\
&\otimes(\beta^{-1}(a_{2})h_{1(1)})(h_{21}\cdot
\beta^{-1}(\beta^{-1}(b_{2})k_{1(1)}))\#\alpha(h_{22})k_{2}\\
=&\sum a_{1}(\alpha(h_{1(0)1})\cdot\beta^{-1}(b_{1}))\#\alpha^{2}(h_{1(0)2})\alpha(k_{1(0)})\\
&\otimes(\beta^{-1}(a_{2})h_{1(1)})(h_{21}\cdot
\beta^{-2}(b_{2})\beta^{-1}(k_{1(1)}))\#\alpha(h_{22})k_{2},
\end{align*}
we have
\begin{align}
&\sum a_{1}(h_{1}\cdot\beta^{-1}(b))_{1}\#\alpha((\alpha(h_{21})k_{1})_{(0)})\otimes\beta^{-1}(a_{2}(h_{1}\cdot\beta^{-1}(b))_{2})(\alpha(h_{21})k_{1})_{(1)}\#\alpha(h_{22})k_{2}\nonumber\\
=&\sum a_{1}(\alpha(h_{1(0)1})\cdot\beta^{-1}(b_{1}))\#\alpha^{2}(h_{1(0)2})\alpha(k_{1(0)})\nonumber\\
&\otimes(\beta^{-1}(a_{2})h_{1(1)})(h_{21}\cdot
\beta^{-2}(b_{2})\beta^{-1}(k_{1(1)}))\#\alpha(h_{22})k_{2}.
\end{align}

Let $a=1$ and $k=1$ in (2.5). Then
\begin{align}
&\sum\beta((h_{1}\cdot\beta^{-1}(b))_{1})\#\alpha((\alpha^{2}(h_{21}))_{(0)})\otimes(h_{1}\cdot\beta^{-1}(b))_{2}\alpha^{2}(h_{21})_{(1)}\#\alpha^{2}(h_{22})\nonumber\\
=&\sum\alpha^{2}(h_{1(0)1})\cdot b_{1}\#\alpha^{3}(h_{1(0)2})\otimes\beta(h_{1(1)})(h_{21}\cdot
\beta^{-1}(b_{2})\#\alpha^{2}(h_{22}).
\end{align}

Applying $id\otimes\varepsilon_{H}\otimes id\otimes\varepsilon_{H}$ to both sides of (2.6), we have
$$\begin{aligned}
&\sum\beta((h_{1}\cdot\beta^{-1}(b))_{1})\varepsilon_{H}(h_{21})\otimes(h_{1}\cdot\beta^{-1}(b))_{2}1_{H}\varepsilon_{H}(h_{22})\\
=&\sum\alpha^{2}(h_{1(0)1})\cdot b_{1}\varepsilon_{H}(h_{1(0)2})\otimes\beta(h_{1(1)})(h_{21}\cdot
\beta^{-1}(b_{2}))\varepsilon_{H}(h_{22}),
\end{aligned}$$

so
$$\sum(h\cdot b)_{1}\otimes(h\cdot b)_{2}=\sum\alpha(h_{1(0)}\cdot b_{1})\otimes\beta(h_{1(1)})(\alpha^{-1}(h_{2})\cdot
\beta^{-1}(b_{2})).$$
Applying $\varepsilon_{B}\otimes id\otimes id\otimes\varepsilon_{H}$ to both sides of (2.6), we have
$$
\sum\alpha^{2}(h_{2(0)})\otimes\beta^{-1}(h_{1}\cdot\beta^{-1}(b))\beta(h_{2(1)})
=\sum\alpha^{2}(h_{1(0)})\otimes\beta(h_{1(1)})\beta^{-1}(h_{2}\cdot
\beta^{-1}(b)),
$$

so
$$\sum h_{2(0)}\otimes (h_{1}\cdot b)\beta^{2}(h_{2(1)})
=\sum h_{1(0)}\otimes\beta^{2}(h_{1(1)})(h_{2}\cdot b).
$$

Let $a=b=1$  in (2.5). Then
$$\begin{aligned}
&\sum1\#\varepsilon_{H}(h_{1})\alpha((\alpha(h_{21})k_{1})_{(0)})\otimes\beta(\alpha(h_{21})k_{1})_{(1)}\#\alpha(h_{22})k_{2}\\
=&\sum1\#\varepsilon_{H}(h_{1(0)1})\alpha^{2}(h_{1(0)2})\alpha(k_{1(0)})\otimes\beta(h_{1(1)})(h_{21}\cdot
k_{1(1)})\#\alpha(h_{22})k_{2},
\end{aligned}$$

and more,
\begin{align}
&\sum1\#\alpha((h_{1}k_{1})_{(0)})\otimes\beta((h_{1}k_{1})_{(1)})\#\alpha(h_{2})k_{2}\nonumber\\
=&\sum1\#\alpha(h_{1(0)})\alpha(k_{1(0)})\otimes\beta(h_{1(1)})(h_{21}\cdot
k_{1(1)})\#\alpha(h_{22})k_{2}.
\end{align}
Applying $\varepsilon_{B}\otimes id\otimes id\otimes\varepsilon_{H}$ to both sides of the above equation, we have
$$\sum(hk)_{(0)}\otimes(hk)_{(1)}=\sum\alpha(h_{1(0)})k_{(0)}\otimes\beta(h_{1(1)})(\alpha^{-1}(h_{2})\cdot \beta^{-1}(k_{(1)})).$$
Conversely if the condition (1), (2), (3) (4) and (5) are satisfied,
\begin{align*}
&\Delta((a\#h)(b\#k))\\
=&\sum \Delta(a(h_{1}\cdot\beta^{-1}(b)\#\alpha(h_{2})k)\\
=&\sum a_{1}(h_{1}\cdot\beta^{-1}(b))_{1}\#\alpha((\alpha(h_{21})k_{1})_{(0)})\otimes\beta^{-1}(a_{2}(h_{1}\cdot\beta^{-1}(b))_{2})(\alpha(h_{21})k_{1})_{(1)}\#\alpha(h_{22})k_{2}\\
=&\sum a_{1}(\alpha(h_{11(0)})\cdot\beta^{-1}(b_{1}))\#\alpha((\alpha(h_{21})k_{1})_{(0)})\\
&\otimes\beta^{-1}(a_{2}(\beta(h_{1(1)})(\alpha^{-1}(h_{12})\cdot \beta^{-2}(b_{2}))))(\alpha(h_{21})k_{1})_{(1)}\#\alpha(h_{22})k_{2}\ by\ (1)\\
=&\sum a_{1}(\alpha(h_{11(0)})\cdot\beta^{-1}(b_{1}))\#\alpha(\alpha(\alpha(h_{211(0)}))k_{1(0)})\otimes[\beta^{-1}(a_{2})(h_{11(1)}(\alpha^{-2}(h_{12})\cdot \beta^{-3}(b_{2})))]\\
&[\beta(\beta(h_{211(1)}))(\alpha^{-1}(\alpha(h_{212}))\cdot\beta^{-1}(k_{1(1)})]\#\alpha(h_{22})k_{2}\\
=&\sum a_{1}(\alpha(h_{11(0)})\cdot\beta^{-1}(b_{1}))\#\alpha^{3}(h_{211(0)})\alpha(k_{1(0)})\otimes[\beta^{-1}(a_{2})(h_{11(1)}(\alpha^{-2}(h_{12})\cdot \beta^{-3}(b_{2})))]\\
& [\beta^{2}(h_{211(1)})(h_{212}\cdot\beta^{-1}(k_{1(1)}))]\#\alpha(h_{22})k_{2}\ by\ (5)\\
=&\sum a_{1}(\alpha(h_{11(0)})\cdot\beta^{-1}(b_{1}))\#\alpha^{3}(h_{212(0)})\alpha(k_{1(0)})\otimes[\beta^{-1}(a_{2})(h_{11(1)}(\alpha^{-2}(h_{12})\cdot \beta^{-3}(b_{2})))]\\
&[(h_{211}\cdot\beta^{-1}(k_{1(1)}))\beta^{2}(h_{212(1)})]\#\alpha(h_{22})k_{2}\ by\ (4)\\
=&\sum a_{1}(h_{1(0)}\cdot\beta^{-1}(b_{1}))\#\alpha^{3}(h_{212(0)})\alpha(k_{1(0)})\otimes[\beta^{-1}(a_{2})(\beta^{-1}(h_{1(1)})(h_{2111}\cdot \beta^{-3}(b_{2})))]\\
&[(\alpha(h_{2112})\cdot\beta^{-1}(k_{1(1)}))\beta^{2}(h_{212(1)})]\#\alpha(h_{22})k_{2}\ by\ (2.3)\\
=&\sum a_{1}(h_{1(0)}\cdot\beta^{-1}(b_{1}))\#\alpha^{3}(h_{212(0)})\alpha(k_{1(0)})\otimes[(\beta^{-2}(a_{2})\beta^{-1}(h_{1(1)}))(\alpha(h_{2111})\cdot \beta^{-2}(b_{2})))]\\
&[(\alpha(h_{2112})\cdot\beta^{-1}(k_{1(1)}))\beta^{2}(h_{212(1)})]\#\alpha(h_{22})k_{2}\\
=&\sum a_{1}(h_{1(0)}\cdot\beta^{-1}(b_{1}))\#\alpha^{3}(h_{212(0)})\alpha(k_{1(0)})\\
&\otimes[(\beta^{-2}(a_{2})\beta^{-1}(h_{1(1)}))(h_{211}\cdot \beta^{-3}(b_{2})\beta^{-2}(k_{1(1)})))] \beta^{3}(h_{212(1)})\#\alpha(h_{22})k_{2}\\
=&\sum a_{1}(h_{1(0)}\cdot\beta^{-1}(b_{1}))\#\alpha^{3}(h_{212(0)})\alpha(k_{1(0)})\\
&\otimes(\beta^{-1}(a_{2})h_{1(1)})[(h_{211}\cdot \beta^{-3}(b_{2})\beta^{-2}(k_{1(1)}))\beta^{2}(h_{212(1)})]\#\alpha(h_{22})k_{2}\\
=&\sum a_{1}(h_{1(0)}\cdot\beta^{-1}(b_{1}))\#\alpha^{3}(h_{211(0)})\alpha(k_{1(0)})\\
&\otimes(\beta^{-1}(a_{2})h_{1(1)})[\beta^{2}(h_{211(1)})(h_{212}\cdot \beta^{-3}(b_{2})\beta^{-2}(k_{1(1)})]\#\alpha(h_{22})k_{2}\ by\ (4)\\
=&\sum a_{1}(\alpha(h_{11(0)})\cdot\beta^{-1}(b_{1}))\#\alpha^{3}(\alpha^{-1}(h_{12(0)}))\alpha(k_{1(0)})\otimes(\beta^{-1}(a_{2})\beta(h_{11(1)}))[\beta(h_{12(1)})\\
&(\alpha^{-1}(h_{21})\cdot \beta^{-3}(b_{2})\beta^{-2}(k_{1(1)})]
\#\alpha(h_{22})k_{2}\ by\ (2.4)\\ =&\sum a_{1}(\alpha(h_{11(0)})\cdot\beta^{-1}(b_{1}))\#\alpha^{3}(\alpha^{-1}(h_{12(0)}))\alpha(k_{1(0)})\\
&\otimes[\beta^{-1}(a_{2})(h_{11(1)}h_{12(1)})]
(h_{21}\cdot \beta^{-2}(b_{2})\beta^{-1}(k_{1(1)})]\#\alpha(h_{22})k_{2}\\
=&\sum a_{1}(\alpha(h_{1(0)1})\cdot\beta^{-1}(b_{1}))\#\alpha^{2}(h_{1(0)2})\alpha(k_{1(0)})\\
&\otimes(\beta^{-1}(a_{2})h_{1(1)})
(h_{21}\cdot \beta^{-2}(b_{2})\beta^{-1}(k_{1(1)})) \#\alpha(h_{22})k_{2}\ by\ comodule\ coalgebra \\
=&\Delta(a\#h)\Delta(b\#k)
\end{align*}
And it is easy to check that $\varepsilon((a\#h)(b\#k))=\varepsilon(a\#h)\varepsilon(b\#k)$. Therefore $(B\#H,\beta\#\alpha)$ is a monoidal Hom-bialgebra.
The proof is completed. $\hfill \Box$
\\

\noindent{\bf Theorem 2.3.}
Let $(B,\beta)$ and $(H,\alpha)$ be monoidal Hom-Hopf algebras, and we have the bicrossproduct monoidal Hom-bialgebra $(B\#H,\beta\#\alpha)$ defined as above. Define the antipode by
\begin{equation}
S(a\#h)=\sum(1_{B}\#S_{H}(h_{(0)}))(S_{B}(\beta^{-2}(a)\beta^{-1}(h_{(1)}))\#1_{H}),
\end{equation}
then $(B\#H,\beta\#\alpha)$ is a monoidal Hom-Hopf algebra (we will call this a  Majid's bicrossproduct).

\noindent{\bf Proof}~~For any $a\#h\in B\#H$,
$$\begin{aligned}
&\sum S(a_{1}\#\alpha(h_{1(0)}))(\beta^{-1}(a_{2})h_{1(1)}\# h_{2})\\
=&\sum [(1_{B}\#S_{H}\alpha(h_{1(0)(0)}))(S_{B}(\beta^{-2}(a_{1})\beta^{-1}\beta(h_{1(0)(1)}))\#1_{H})](\beta^{-1}(a_{2})h_{1(1)}\# h_{2})\\
=&\sum [(1_{B}\#S_{H}\alpha(h_{1(0)(0)}))(S_{B}(\beta^{-2}(a_{1})h_{1(0)(1)})\#1_{H})](\beta^{-1}(a_{2})h_{1(1)}\otimes h_{2})\\
=&\sum (1_{B}\#S_{H}\alpha^{2}(h_{1(0)(0)}))[(S_{B}(\beta^{-2}(a_{1})h_{1(0)(1)})\#1_{H})(\beta^{-2}(a_{2})\beta^{-1}(h_{1(1)})\#\alpha^{-1}(h_{2}))]\\
=&\sum (1_{B}\#S_{H}\alpha^{2}(h_{1(0)(0)}))[(S_{B}(\beta^{-2}(a_{1})h_{1(0)(1)})(\beta^{-2}(a_{2})\beta^{-1}(h_{1(1)})\# h_{2}]\\
=&\sum (1_{B}\#S_{H}\alpha^{2}(h_{1(0)(0)}))[(S_{B}(h_{1(0)(1)})S_{B}\beta^{-2}(a_{1}))(\beta^{-2}(a_{2})\beta^{-1}(h_{1(1)})\# h_{2}]\\
=&\sum (1_{B}\#S_{H}\alpha^{2}(h_{1(0)(0)}))[(S_{B}(h_{1(0)(1)})(S_{B}\beta^{-3}(a_{1})\beta^{-3}(a_{2})))h_{1(1)}\# h_{2}]\\
=&\sum (1_{B}\#S_{H}\alpha^{2}(h_{1(0)(0)}))[(S_{B}\beta(h_{1(0)(1)})h_{1(1)}\# h_{2}]\varepsilon_{B}(a)\\
=&\sum (1_{B}\#S_{H}\alpha(h_{1(0)}))[(S_{B}\beta(h_{1(1)1})\beta(h_{1(1)2})\# h_{2}]\varepsilon_{B}(a)\\
=&\sum (1_{B}\#S_{H}\alpha(h_{1(0)}))(\varepsilon_{B}(h_{1(1)}))1\# h_{2})\varepsilon_{B}(a)\\
=&\sum (1_{B}\#S_{H}(h_{1})h_{2})\varepsilon_{B}(a)\\
=&\sum (1_{B}\#1_{H})\varepsilon_{B}(a)\varepsilon(h).
\end{aligned}$$
That is $S\ast id=\varepsilon$. Similarly we have $id\ast S=\varepsilon$.
This completes the proof. $\hfill \Box$
\\

\noindent{\bf Remark 2.4.}
In Theorem 2.2, if the module action is trivial, that is $h\cdot a=\varepsilon_{H}(h)\alpha(a)$, we have a bicrossproduct with monoidal Hom-smash coproduct and tensor product. If the comodule is trivial, that is $\rho(h)=\alpha^{-1}(h)\otimes 1$, we have a bicrossproduct with monoidal  Hom-smash product and tensor coproduct.
\\

\noindent{\bf Example 2.5.}
Let $B=span\{1,x\}$ over a fixed field $k$ with $chark\neq2$. Then define $\beta$ as a $k$-linear automorphism of $B$ by
$$\beta(1_{B})=1_{B},\ \ \beta(x)=-x.$$
Define the multiplication on $B$ by
$$1_{B}1_{B}=1_{B},\ 1_{B}x=-x,\ x^{2}=0,$$
then it is not hard to check that  $(B,\beta)$ is a Hom-associative algebra.

 For $B$,  define the comultiplication and the counit by
\begin{align*}
&\Delta(1_{B})=1_{B}\otimes 1_{B},\ \Delta(x)=(-x)\otimes 1+1\otimes(-x),\\
&\varepsilon(1_{B})=1,\ \varepsilon(x)=0,\\
&S_{B}(1_{B})=1_{B},\ S_{B}(x)=-x,
\end{align*}
then it is easy to check that $(B,\beta)$ is a monoidal Hom-Hopf algebra.

Let $H=span\{1,g|g^{2}=1\}$ be the group algebra. Obviously $(H,id)$ is a monoidal Hom-Hopf algebra.

Define the left action of $H$ on $B$ by $\cdot:H\otimes B\rightarrow B$ by
$$1_{H}\cdot 1_{B}=1_{B},\ 1_{H}\cdot x=-x,\ g\cdot 1_{B}=1_{B},\ g\cdot x=x.$$
It is easy to check that $(B,\beta)$ is a left $(H,id)$-module algebra.

Define the right coaction of $B$ on $H$ $\rho:H\rightarrow H\otimes B$ by
$$\rho(1_{H})=1_{H}\otimes 1_{B},\ \rho(g)=g\otimes 1_{B},$$
We could verify  that $(H,id)$ is a right $(B,\beta)$-comodule coalgebra, and the conditions in the above proposition are satisfied. Hence we have a monoidal Hom-Hopf algebra $(B\#H,\beta\#id)$ (the Majid's bicrossproduct), where
\begin{align*}
&\Delta(1_{B}\#1_{H})=1_{B}\#1_{H}\otimes 1_{B}\#1_{H},\ \varepsilon(1_{B}\#1_{H})=1,\\
&\Delta(1_{B}\#g)=1_{B}\#g\otimes1_{B}\#g,\ \varepsilon(1_{B}\#g)=1,\\
&\Delta(x\#1_{H})=(-x)\# 1_{H}\otimes1_{B}\#1_{H},\ \varepsilon(x\#1_{H})=0,\\
&\Delta(x\#g)=(-x)\#g\otimes1_{B}\#g+1_{B}\#g\otimes(-x)\#g,\ \varepsilon(x\#g)=0.\\
&S(1_{B}\#1_{H})=1_{B}\#1_{H},\ S(1_{B}\#g)=1_{B}\#g,\ S(x\#1_{H})=x\otimes 1_{H},\ S(x\#g)=(-x)\#g.
\end{align*}
As to the multiplication, it is the same as in the Example 3.7 in \cite{LS}.

\section{A class of monoidal Hom-Hopf algebras}
\def\theequation{3.\arabic{equation}}
\setcounter{equation} {0}

In this subsection, we will construct a class of monoidal Hom-Hopf algebras as an application of the theory constructed
 in Section 2.

\noindent{\bf Lemma 3.1.} Let $(H,\alpha)$ be a monoidal Hom-Hopf algebra.  Define the left action of $(H^{op},\alpha)$ on $(H,\alpha)$ by
$$\cdot:H^{op}\otimes H\rightarrow H,\ h\cdot a=\sum(S(h_{1})\alpha^{-1}(a))\alpha(h_{2}),$$
for any $h,a\in H$. Then $(H,\alpha)$ is a left $(H^{op},\alpha)$-module algebra.

\noindent{\bf Proof}~~For any $h,k,a,b\in H$, firstly, $1\cdot a=\alpha(a)$, and
$$\begin{aligned}
\alpha(h\cdot a)=\sum (S\alpha(h_{1})a)\alpha^{2}(h_{2})=\alpha(h)\cdot\alpha( a).
\end{aligned}$$
Secondly
$$\begin{aligned}
(hk)\cdot a&=\sum(S(h_{1}k_{1})\alpha^{-1}(a))\alpha(h_{2}k_{2})\\
           &=\sum[(S(k_{1})S(h_{1}))\alpha^{-1}(a)]\alpha(h_{2}k_{2})\\
           &=\sum[\alpha S(k_{1})(S(h_{1})\alpha^{-2}(a))]\alpha(h_{2}k_{2})\\
           &=\sum[(S(k_{1})(\alpha^{-1}S(h_{1})\alpha^{-3}(a)))\alpha(h_{2})]\alpha^{2}(k_{2})\\
           &=\sum[\alpha S(k_{1})((\alpha^{-1}S(h_{1})\alpha^{-3}(a))h_{2})]\alpha^{2}(k_{2})\\
           &=\sum\alpha(k)\cdot[(S(h_{1})\alpha^{-2}(a))\alpha(h_{2})]\\
           &=\alpha(k)\cdot[h\cdot\alpha^{-1}(a)],
\end{aligned}$$
Now $(H,\alpha)$ is a left $(H^{op},\alpha)$-module. Finally

$$\begin{aligned}
\sum (h_{1}\cdot a)(h_{2}\cdot b)&=\sum ((S(h_{11})\alpha^{-1}(a))\alpha(h_{12}))((S(h_{21})\alpha^{-1}(b))\alpha(h_{22}))\\
                                 &=\sum((S(h_{11})\alpha^{-1}(a))\alpha(h_{12}))(S\alpha(h_{21})(\alpha^{-1}(b)h_{22}))\\
                                 &=\sum[(S(h_{11})\alpha^{-1}(a))(h_{12}S(h_{21}))](b\alpha(h_{22}))\\
                                 &=\sum[(S\alpha^{-1}(h_{1})\alpha^{-1}(a))(\alpha(h_{211})S\alpha(h_{212}))](b\alpha(h_{22}))\\
                                 &=\sum[(S(h_{1})a)](bh_{2})=(S(h_{1})\alpha^{-1}(ab))\alpha(h_{2})\\
                                 &=h\cdot(ab),
\end{aligned}$$

and
$$h\cdot 1=\sum\alpha(S(h_{1}))\alpha(h_{2})=\varepsilon(h)1.$$
This completes the proof.$\hfill \Box$
\\

\noindent{\bf Lemma 3.2.}
Let $(H,\alpha)$ be a monoidal Hom-Hopf algebra. Define the right coaction of $(H,\alpha)$ on $(H^{op},\alpha)$ by
$$\rho:H^{op}\rightarrow H^{op}\otimes H,\ \rho(h)=\sum\alpha(h_{12})\otimes S(h_{11})\alpha^{-1}(h_{2}),$$
for any $h\in H$. Then $(H^{op},\alpha)$ is a right $(H,\alpha)$-comodule coalgebra.

\noindent{\bf Proof}~~For any $h\in H$, firstly
$$\sum\alpha(h_{12})\varepsilon(h_{11})\varepsilon(h_{2})=\alpha^{-1}(h),$$
and
$$\rho(\alpha(h))=\sum\alpha(\alpha(h)_{12})\otimes S(\alpha(h)_{11})\alpha^{-1}(\alpha(h)_{2})=(\alpha\otimes\alpha)\rho(h).$$

Secondly
$$\begin{aligned}
&\sum\alpha^{-1}(\alpha(h_{12}))\otimes\Delta(S(h_{11})\alpha^{-1}(h_{2}))\\
=&\sum h_{12}\otimes S(h_{112})\alpha^{-1}(h_{21})\otimes S(h_{111})\alpha^{-1}(h_{22})\\
=&\sum \alpha^{2}(h_{1122})\otimes S\alpha(h_{1121})\alpha^{-1}(h_{12})\otimes S(h_{111})\alpha^{-2}(h_{2})\\
=&\sum \alpha^{2}(h_{1212})\otimes S\alpha(h_{1211})h_{122}\otimes S\alpha^{-1}(h_{11})\alpha^{-2}(h_{2})\\
=&\sum h_{(0)(0)}\otimes h_{(0)(1)}\otimes\alpha^{-1}(h_{(1)}).
\end{aligned}$$
Then $(H^{op},\alpha)$ is a right $(H,\alpha)$-comodule. Finally

$$\begin{aligned}
&\sum h_{1(0)}\otimes h_{2(0)}\otimes h_{1(1)}h_{1(1)}\\
=&\sum \alpha(h_{112})\otimes\alpha(h_{212})\otimes [S(h_{111})\alpha^{-1}(h_{12})][S(h_{211})\alpha^{-1}(h_{22})]\\
=&\sum \alpha(h_{112})\otimes\alpha(h_{212})\otimes [S(h_{111})(\alpha^{-2}(h_{12})\alpha^{-1}S(h_{211}))]h_{22}\\
=&\sum h_{12}\otimes\alpha(h_{212})\otimes [S\alpha^{-1}(h_{11})(h_{2111}S(h_{2112}))]h_{22}\\
=&\sum h_{12}\otimes h_{21}\otimes S(h_{11})h_{22}\\
=&\sum\alpha(h_{121})\otimes\alpha(h_{122})\otimes S(h_{11})\alpha^{-1}(h_{2})\\
=&\sum\Delta(\alpha(h_{12}))\otimes S(h_{11})\alpha^{-1}(h_{2}),
\end{aligned}$$
and
$$\sum \varepsilon(\alpha(h_{12}))S(h_{11})\alpha^{-1}(h_{2})=\varepsilon(h)1.$$
This completes the proof. $\hfill \Box$
\\

Now we have the following result.

\noindent{\bf Theorem 3.3.}
Let $(H,\alpha)$ be a monoidal Hom-Hopf algebra. The action $\cdot:H^{op}\otimes H\rightarrow H$ and coaction $\rho:H^{op}\rightarrow H^{op}\otimes H$ are defined as above.
 Then we have  Majid's bicrossproduct  monoidal Hom-Hopf algebra $(H\ast H^{op},\alpha\ast\alpha)$.

\noindent{\bf Proof}~~For any $h,k,a,b\in H$, firstly
$$\begin{aligned}
\Delta(h\cdot a)&=\sum\Delta((S(h_{1})\alpha^{-1}(a))\alpha(h_{2}))\\
                &=\sum (S(h_{12})\alpha^{-1}(a_{1}))\alpha(h_{21})\otimes(S(h_{11})\alpha^{-1}(a_{2}))\alpha(h_{22})\\
                &=\sum \alpha(h_{21})\cdot a_{1}\otimes(S\alpha^{-1}(h_{1})\alpha^{-1}(a_{2}))\alpha(h_{22})\\
                &=\sum \alpha(h_{12})\cdot a_{1}\otimes(S(h_{11})\alpha^{-1}(a_{2}))h_{2}\\
                &=\sum \alpha(h_{12})\cdot a_{1}\otimes S\alpha(h_{11})(\alpha^{-1}(a_{2})\alpha^{-1}(h_{2}))\\
                &=\sum \alpha(h_{12})\cdot a_{1}\otimes S\alpha(h_{11})(\varepsilon(h_{21})\alpha^{-1}(a_{2})h_{22})\\
                &=\sum \alpha^{2}(h_{112})\cdot a_{1}\otimes \alpha^{2}S(h_{111})(((\alpha^{-2}(h_{12})S\alpha^{-2}(h_{21}))\alpha^{-2}(a_{2}))h_{22})\\
                &=\sum \alpha^{2}(h_{112})\cdot a_{1}\otimes \alpha^{2}S(h_{111})((\alpha^{-1}(h_{12})(S\alpha^{-2}(h_{21})\alpha^{-3}(a_{2})))h_{22})\\
                &=\sum \alpha^{2}(h_{112})\cdot a_{1}\otimes \alpha^{2}S(h_{111})(h_{12}((S\alpha^{-2}(h_{21})\alpha^{-3}(a_{2}))\alpha^{-1}(h_{22})))\\
                &=\sum \alpha^{2}(h_{112})\cdot a_{1}\otimes \alpha^{2}S(h_{111})(h_{12}(\alpha^{-2}(h_{2})\cdot\alpha^{-2}(a_{2})))\\
                &=\sum \alpha^{2}(h_{112})\cdot a_{1}\otimes \alpha(S(h_{111})\alpha^{-1}(h_{12}))(\alpha^{-1}(h_{2})\cdot\alpha^{-1}(a_{2}))\\
                &=\sum \alpha(h_{1(0)})\cdot a_{1}\otimes\alpha(h_{1(1)})(\alpha^{-1}(h_{2})\cdot\alpha^{-1}(a_{2})),
\end{aligned}$$
and $$\varepsilon(h\cdot a)=\sum \varepsilon(h_{1})\varepsilon(h_{2})\varepsilon(a)=\varepsilon(h)\varepsilon(a).$$

Secondly
$$\begin{aligned}
&\sum\alpha(h_{212})\otimes(h_{1}\cdot b)\alpha^{2}(S(h_{211})\alpha^{-1}(h_{22}))\\
=&\sum\alpha(h_{212})\otimes((S(h_{11})\alpha^{-1}(b))\alpha(h_{12}))\alpha^{2}(S(h_{211})\alpha^{-1}(h_{22}))\\
=&\sum\alpha(h_{212})\otimes[(S(h_{11})\alpha^{-1}(b))(h_{12} S\alpha(h_{211}))]\alpha^{2}(h_{22})\\
=&\sum\alpha(h_{212})\otimes[(S(h_{1})b)\varepsilon(h_{211})]\alpha^{2}(h_{22})\\
=&\sum h_{21}\otimes[(S(h_{1})b)]\alpha^{2}(h_{22})\\
=&\sum h_{12}\otimes S\alpha^{2}(h_{11})(bh_{2})\\
=&\sum h_{12}\otimes S\alpha^{2}(h_{11})\varepsilon(h_{21})(b\alpha(h_{22}))\\
=&\sum \alpha(h_{112})\otimes (S\alpha^{2}(h_{111})(h_{12}S(h_{21})))(b\alpha(h_{22}))\\
=&\sum \alpha(h_{112})\otimes [\alpha^{2}S(h_{111})\alpha(h_{12})][\alpha S(h_{21})(\alpha^{-1}(b)h_{22})]\\
=&\sum \alpha(h_{112})\otimes \alpha^{2}(S(h_{111})\alpha^{-1}(h_{12}))(h_{2}\cdot b)\\
=&\sum h_{1(0)}\otimes\beta^{2}(h_{1(1)})(h_{2}\cdot b),
\end{aligned}$$

Finally
$$\begin{aligned}
&\sum (h\circ k)_{(0)}\otimes (h\circ k)_{(1)}\\
=&\sum\alpha(k_{12}h_{12})\otimes[S(h_{11})S(k_{11})][\alpha^{-1}(k_{2})\alpha^{-1}(h_{2})] \\
=&\sum\alpha(k_{12})\alpha(h_{12})\otimes S\alpha(h_{11})[(S\alpha^{-1}(k_{11})\alpha^{-2}(k_{2}))\alpha^{-1}(h_{2})]\\
=&\sum\alpha(k_{12})\alpha^{2}(h_{112})\otimes [S\alpha(h_{111})(\alpha^{-1}(h_{12})\alpha^{-1}S(h_{21}))][(S\alpha^{-1}(k_{11})\alpha^{-2}(k_{2}))h_{22}]\\
=&\sum\alpha(k_{12})\alpha^{2}(h_{112})\otimes
[S\alpha(h_{111})h_{12}][S(h_{21})((S\alpha^{-2}(k_{11})\alpha^{-3}(k_{2}))\alpha^{-1}(h_{22}))]\\
=&\sum\alpha(k_{12})\alpha^{2}(h_{112})\otimes
[S\alpha(h_{111})h_{12}][(S\alpha^{-1}(h_{21})(S\alpha^{-2}(k_{11})\alpha^{-3}(k_{2})))h_{22}]\\
=&\sum\alpha(k_{12})\alpha^{2}(h_{112})\otimes
\alpha[S(h_{111})\alpha^{-1}(h_{12})][\alpha^{-1}(h_{2})\cdot\alpha^{-1}(S(k_{11})\alpha^{-1}(k_{2}))]\\
=&\sum\alpha(h_{1(0)})\circ k_{(0)}\otimes\alpha(h_{1(1)})(\alpha^{-1}(h_{2})\cdot\alpha^{-1}(k_{(1)})).
\end{aligned}$$
This completes the proof. $\hfill \Box$
\\

By the above proposition, the Hopf algebra structure on $(H\ast H^{op},\alpha\ast\alpha)$ is  given by
\begin{align*}
&(a\#h)(b\#k)=\sum a[(S(h_{11})\beta^{-2}(b))\alpha(h_{12})]\#k\alpha(h_{2}),\\
&\Delta(a\#h)=\sum a_{1}\#\alpha^{2}(h_{112})\otimes\beta^{-1}(a_{2})(S(h_{111})\alpha^{-1}(h_{12}))\#h_{2},
\end{align*}
for any $a,b,h,k\in H.$

\section{Drinfel'd Double of monoidal Hom-Hopf algebras}
\def\theequation{4.\arabic{equation}}
\setcounter{equation} {0}

In this section, we will introduce Drinfel'd double of monoidal Hom-Hopf algebras. The monoidal Hom-Hopf algebras in this section are all assumed to be finite-dimensional. First we have the following definition:

\noindent{\bf Definition 4.1.}
Let $(B,\beta)$ and $(H,\alpha)$ be two monoidal Hom-Hopf algebras. Then $((B,\beta),(H,\alpha))$ is called
 a {\sl matched monoidal Hom-Hopf algebra pair} if $H$ is a right $B$-Hom-module coalgebra and $B$ is a left $H$-Hom-module coalgebra via
$$\triangleleft:H\otimes B\rightarrow H,\ \triangleright:H\otimes B\rightarrow B,$$
with the compatible conditions:
\begin{enumerate}
\item
[(1)]$(hg)\triangleleft a=\sum(h\triangleleft(g_{1}\triangleright\beta^{-1}(a_{1})))(\alpha(g_{2})\triangleleft a_{2}),\ 1_{H}\triangleleft a=\varepsilon_{B}(a)1_{H}$,
\item
[(2)]$h\triangleright(ab)=\sum(h_{1}\triangleright\beta(a_{1}))((\alpha^{-1}(h_{2})\triangleleft a_{2})\triangleright b),\ 1_{B}\triangleleft h=\varepsilon_{H}(h)1_{B}$,
\item
[(3)]$\sum(h_{1}\triangleleft a_{1})\otimes(h_{2}\triangleright a_{2})=(h_{2}\triangleleft a_{2})\otimes(h_{1}\triangleright a_{1})$.
\end{enumerate}

By the above definition, we have

\noindent{\bf Theorem 4.2.}
Given a matched monoidal Hom-Hopf algebra pair $((B,\beta),(H,\alpha))$, we have a double cross product monoidal Hom-Hopf algebra $(B\bowtie H,\beta\otimes\alpha)$ built on the vector space $B\otimes H$ with the product
$$(a\ltimes h)(b\ltimes g)=\sum a(h_{1}\triangleright b_{1})\ltimes(h_{2}\triangleleft b_{2})g$$
and the usual tensor coproduct.

The antipode is defined by
$$S(a\ltimes h)=\sum S_{H}(h_{2})\triangleright S_{B}(a_{2})\ltimes S_{H}(h_{1})\triangleleft S_{B}(a_{1}).$$

\noindent{\bf Proof}~~Firstly the unit and counit are easy to verify.
Secondly
$$\begin{aligned}
&(\beta(a)\ltimes\alpha(h))((b\ltimes g)(c\ltimes f))\\
=&\sum(\beta(a)\ltimes\alpha(h))(b(g_{1}\triangleright c_{1})\ltimes(g_{2}\triangleleft c_{2})f)\\
=&\sum \beta(a)(\alpha(h_{1})\triangleright(b_{1}(g_{11}\triangleright c_{11})))\ltimes(\alpha(h_{2})\triangleleft (b_{2}(g_{12}\triangleright c_{12})))((g_{2}\triangleleft c_{2})f)\\
=&\sum \beta(a)[(\alpha(h_{11})\triangleright\beta(b_{11}))((h_{12}\triangleleft b_{12})\triangleright(g_{11}\triangleright c_{11})))]\\
&\ltimes((h_{2}\triangleleft b_{2})\triangleleft(\alpha(g_{12})\triangleright \beta(c_{12})))((g_{2}\triangleleft c_{2})f)\ (b)\\
=&\sum \beta(a)[(\alpha(h_{11})\triangleright\beta(b_{11}))((h_{12}\triangleleft b_{12})\triangleright(\alpha^{-1}(g_{1})\triangleright \beta^{-1}(c_{1}))))]\\
&\ltimes((h_{2}\triangleleft b_{2})\triangleleft(\alpha(g_{21})\triangleright \beta(c_{21})))((\alpha(g_{22})\triangleleft \beta(c_{22}))f)\\
\end{aligned}$$
$$\begin{aligned}
=&\sum \beta(a)[(h_{1}\triangleright b_{1})(((\alpha^{-1}(h_{21})\triangleleft \beta^{-1}(b_{21}))\alpha^{-1}(g_{1}))\triangleright c_{1})]\\
&\ltimes[((h_{22}\triangleleft b_{22})\triangleleft(g_{21}\triangleright c_{21}))(\alpha(g_{22})\triangleleft \beta(c_{22}))]\alpha(f)\\
=&\sum \beta(a)[(h_{1}\triangleright b_{1})(((\alpha^{-1}(h_{21})\triangleleft \beta^{-1}(b_{21}))\alpha^{-1}(g_{1}))\triangleright c_{1})]\\
&\ltimes[(h_{22}\triangleleft b_{22})g_{2}\triangleleft \beta(c_{2})]\alpha(f)\ (a)\\
=&\sum (a(h_{1}\triangleright b_{1}))((h_{21}\triangleleft b_{21})g_{1}\triangleright \beta(c_{1}))\ltimes[(h_{22}\triangleleft b_{22})g_{2}\triangleleft \beta(c_{2})]\alpha(f)\ (a)\\
=&\sum (a(h_{1}\triangleright b_{1})\ltimes(h_{2}\triangleleft b_{2})g)(\beta(c)\ltimes\alpha(f))\\
=&((a\ltimes h)(b\ltimes g))(\beta(c)\ltimes\alpha(f)).
\end{aligned}$$
And
\begin{align*}
&\Delta(a\ltimes h)\Delta(b\ltimes g)\\
=&\sum((a_{1}\ltimes h_{1})\otimes(a_{2}\ltimes h_{2}))((b_{1}\ltimes g_{1})\otimes(b_{2}\ltimes g_{2}))\\
=&\sum a_{1}(h_{11}\triangleright b_{11})\ltimes(h_{12}\triangleleft b_{12})g_{1}\otimes a_{2}(h_{21}\triangleright b_{21})\ltimes(h_{22}\triangleleft b_{22})g_{2}\\
=&\sum a_{1}(\alpha^{-1}(h_{1})\triangleright \beta^{-1}(b_{1}))\ltimes(\alpha(h_{211})\triangleleft \beta(b_{211}))g_{1}\\
&\otimes a_{2}(\alpha(h_{212})\triangleright \beta(b_{212}))\ltimes(h_{22}\triangleleft b_{22})g_{2}\\
=&\sum a_{1}(\alpha^{-1}(h_{1})\triangleright \beta^{-1}(b_{1}))\ltimes(\alpha(h_{212})\triangleleft \beta(b_{212}))g_{1}\\
&\otimes a_{2}(\alpha(h_{211})\triangleright \beta(b_{211}))\ltimes(h_{22}\triangleleft b_{22})g_{2}\ (c)\\
=&\sum a_{1}(h_{11}\triangleright b_{11})\ltimes(h_{21}\triangleleft b_{21})g_{1}\otimes a_{2}(h_{12}\triangleright b_{12})\ltimes(h_{22}\triangleleft b_{22})g_{2}\\
=&\Delta((a\ltimes h)(b\ltimes g)).
\end{align*}
Obviously $\varepsilon=\varepsilon_{B}\otimes\varepsilon_{H}$ is a monoidal Hom-algebra map.

It is now straightforward to check that
\begin{align*}
&(\beta\otimes\alpha)((a\ltimes h)(b\ltimes g))=(\beta\otimes\alpha)(a\ltimes h)(\beta\otimes\alpha)(b\ltimes g),\\
&\Delta\circ(\beta\otimes\alpha)=[(\beta\otimes\alpha)\otimes(\beta\otimes\alpha)]\circ\Delta.
\end{align*}
 Thus $B\bowtie H$ is a monoidal Hom-bialgebra. Finally, it is easy to verify that $S\ast id=id\ast S=(1_{B}\ltimes 1_{H})\varepsilon .$

This completes the proof. $\hfill \Box$
\\

\noindent{\bf Proposition 4.3.}
Suppose the bicrossproduct $(B\#H,\beta\otimes\alpha)$ are defined as before. If we define $\triangleright:B^{\ast}\otimes H\rightarrow H$ by
$$f\triangleright h=\sum\langle f,\beta(h_{(1)})\rangle\alpha^{2}(h_{(0)}),$$
and
$\triangleleft:B^{\ast}\otimes H\rightarrow B^{\ast}$ by
$$\langle f\triangleleft h,a\rangle=\langle f,\alpha^{-1}(h)\cdot\beta^{-2}(a)\rangle.$$
Then $(B^{\ast},(\beta^{-1})^{\ast})$ is a right $(H,\alpha)$-Hom module coalgebra, and $(H,\alpha)$ is a left $(B^{\ast},(\beta^{-1})^{\ast})$-Hom module coalgebra.

\noindent{\bf Proof}~~The proof is straightforward. $\hfill \Box$
\\

\noindent{\bf Proposition 4.4.}
Under the action defined in the above proposition, $(H\bowtie B^{\ast},(\alpha\otimes\beta^{-1})^{\ast})$ is a matched monoidal Hom-Hopf algebra pair. The product and antipode are given by
\begin{align*}
&(h\ltimes f)(k\ltimes g)=\sum h\alpha^{2}(k_{1(0)})\ltimes \langle f,\beta(k_{1(1)})(\alpha^{-1}(k_{2})\cdot\beta^{-2}(?))\rangle g,\\
&S(h\ltimes f)=\sum S^{\ast}_{B}(f_{2})\triangleright S_{H}(h_{2})\ltimes S^{\ast}_{B}(f_{1})\triangleleft S_{H}(h_{1}),
\end{align*}
for any $h,k\in H$ and $f,g\in B^{\ast}$.

\noindent{\bf Proof}~~We need only to verify the compatible conditions. Firstly for any $a,b\in B,\ h,k\in H$, and $f,g\in B^{\ast}$,
$$\begin{aligned}
&\sum \langle f\triangleleft(g_{1}\triangleright\alpha^{-1}(h_{1})),a_{1}\rangle\langle (\beta^{-1})^{\ast}(g_{2})\triangleleft h_{2},a_{2}\rangle\\
=&\sum \langle f,\alpha^{-1}((g_{1}\triangleright\alpha^{-1}(h_{1})))\cdot \beta^{-2}(a_{1})\rangle\langle (\beta^{-1})^{\ast}(g_{2})\triangleleft h_{2},a_{2}\rangle\\
=&\sum \langle f,\alpha^{-1}((g_{1}\triangleright\alpha^{-1}(h_{1})))\cdot \beta^{-2}(a_{1})\rangle\langle (\beta^{-1})^{\ast}(g_{2}),\alpha^{-1}(h_{2})\cdot \beta^{-2}(a_{2})\rangle\\
=&\sum \langle f,(\beta^{\ast}(g_{1})\triangleright\alpha^{-2}(h_{1}))\cdot \beta^{-2}(a_{1})\rangle\langle g_{2},\alpha^{-2}(h_{2})\cdot \beta^{-3}(a_{2})\rangle\\
=&\sum \langle f,h_{1(0)}\cdot \beta^{-2}(a_{1})\rangle\langle g_{1},h_{1(1)}\rangle\langle g_{2},\alpha^{-2}(h_{2})\cdot \beta^{-3}(a_{2})\rangle\\
=&\sum \langle f,h_{1(0)}\cdot \beta^{-2}(a_{1})\rangle\langle g,h_{1(1)}(\alpha^{-2}(h_{2})\cdot \beta^{-3}(a_{2}))\rangle\\
=&\sum\langle f,(\alpha^{-1}(h)\cdot\beta^{-2}(a))_{1}\rangle\langle g,(\alpha^{-1}(h)\cdot\beta^{-2}(a))_{2}\rangle\\
=&\langle (f\ast g)\triangleleft h,a\rangle.
\end{aligned}$$

Then

$$\begin{aligned}
f\triangleright(hk)&=\sum\langle f,\beta((hk)_{(1)})\rangle\alpha^{2}((hk)_{(0)})\\
                   &=\sum\langle f,\beta^{2}(h_{1(1)})(g_{2}\cdot k_{(1)})\rangle\alpha^{3}(h_{1(0)})\alpha^{2}(k_{(0)})\\
                   &=\sum\langle f_{1},\beta^{2}(h_{1(1)})\rangle\langle f_{2},h_{2}\cdot k_{(1)}\rangle\alpha^{3}(h_{1(0)})\alpha^{2}(k_{(0)})\\
                   &=\sum (f_{1}\triangleright \beta(h))((\beta^{\ast}(f_{2})\triangleleft h_{2})\triangleright k).
\end{aligned}$$

Finally

$$\begin{aligned}
&\sum \langle f_{1}\triangleleft h_{1},a\rangle\langle g,f_{2}\triangleright h_{2}\rangle\\
=&\sum \langle f_{1},\alpha^{-1}(h_{1})\cdot \beta^{-2}(a)\rangle\langle f_{2},\beta(h_{2(1)})\rangle \langle g,\alpha^{2}(h_{2(0)})\rangle\\
=&\sum \langle f,(\alpha^{-1}(h_{1})\cdot \beta^{-2}(a))\beta(h_{2(1)})\rangle\langle g,\alpha^{2}(h_{2(0)})\rangle\\
=&\sum \langle f_{2},\alpha^{-1}(h_{2})\cdot \beta^{-2}(a)\rangle\langle f_{1},\beta(h_{1(1)})\rangle\langle g,\alpha^{2}(h_{1(0)})\rangle\\
=&\sum \langle f_{2}\triangleleft h_{2},a\rangle\langle g,f_{1}\triangleright h_{1}\rangle.
\end{aligned}$$
This completes the proof.$\hfill \Box$
\\

\noindent{\bf Corollary 4.5.} We have the Drinfel'd double $\mathcal{D}(H)=H^{op}\bowtie H^{\ast}$. Moreover the product is given by
$$(h\ltimes f)(k\ltimes g)=\sum \alpha^{2}(k_{21})h\ltimes \langle f,S(k_{1})(\alpha^{-2}(?)k_{22})\rangle g,$$
for any $h,k\in H$ and $f,g\in H^{\ast}$.

\noindent{\bf Proof}~~For the bicrossproduct monoidal Hom-Hopf algebra $(H\ast H^{op},\alpha\ast\alpha)$ in Theorem 3.3,
 by Proposition 4.3 and Proposition 4.4, it is straightforward to get this result.$\hfill \Box$
 \\

\noindent{\bf Example 4.6.} Let $G$ be a finite group and $\phi$ an automorphism of $G$. Then $(kG,\phi)$ with the following structure map is a monoidal Hom-Hopf algebra:
$$g\cdot h=\phi(gh),\ \Delta(g)=\phi^{-1}(g)\otimes \phi^{-1}(g),\ \varepsilon(g)=1,\ S(g)=g^{-1}.$$

Let $\{e_{g}\}_{g\in G}$ be the dual basis of the basis of $kG$. Then we have the monoidal Hom-Hopf algebra $k^{G}$, dual of $kG$, with the multiplication
$$e_{g}e_{h}=\delta_{g,h}e_{\phi(g)},$$
and the comultiplication, counit and antipode
$$\Delta(e_{g})=\sum_{uv=\phi^{-1}(g)}e_{u}\o e_{v},\ \varepsilon(e_{g})=\delta_{g,1},\ S(e_{g})=e_{g^{-1}},$$
for any $g,h\in G.$

By Corollary 4.5, the multiplication in $D(kG)$ is given by
$$(g\o e_{h})(p\o e_{q})=\phi(pg)\o \delta_{php^{-1},q}\ e_{\phi(q)}.$$

\noindent{\bf Definition 4.7.}
A quasitriangular monoidal Hom-Hopf algebra is a monoidal Hom-Hopf algebra $(H,\alpha)$ with an element $R\in H\otimes H$ satisfying
\begin{enumerate}
\item
[(1)]$(\varepsilon\otimes id)R=(id\otimes\varepsilon)R=1,$
\item
[(2)]$\Delta^{op}(x)R=R\Delta(x)$ for any $x\in H,$
\item
[(3)]$(\Delta\otimes\alpha^{-1})R=R^{13}R^{23},$
\item
[(4)]$(\alpha^{-1}\otimes\Delta)R=R^{13}R^{12}.$
\end{enumerate}

\noindent{\bf Remark 4.8.} Note that the above definition is different from Definition 2.7 in \cite{Y}.
\\

In the above definition, we will denote $R=\sum R^{(1)}\otimes R^{(2)}$ throughout this section. Therefore we could rewrite the equalities (3) and (4)
as follows:
\begin{align}
&(\Delta\otimes id)R=\sum\alpha(R^{(1)})\otimes\alpha(r^{(1)})\otimes\alpha(R^{(2)}r^{(2)}),\\
&(id\otimes\Delta)R=\sum\alpha(R^{(1)}r^{(1)})\otimes\alpha(r^{(2)})\otimes\alpha(R^{(2)}).
\end{align}

\noindent{\bf Example 4.9.} Let $H$ be the monoidal Hom-algebra generated by the elements $1_{H},\ g,\ x$ satisfying the following relations:
\begin{align*}
&1_{H}1_{H}=1_{H},\ 1_{H}g=g1_{H}=g,\ 1_{H}x=x1_{H}=-x,\\
&g^{2}=1_{H},\ x^{2}=0,\ gx=-xg.
\end{align*}
The automorphism $\alpha:H\rightarrow H$ is defined by
$$\alpha(1_{H})=1_{H},\ \alpha(g)=g,\ \alpha(x)=-x,\ \alpha(gx)=-gx.$$
Then $(H,\alpha)$ is a Hom-associative algebra, and $\alpha^{2}=id$.

Define
\begin{align*}
&\Delta(1_{H})=1_{H}\otimes 1_{H},\ \Delta(g)=g\otimes g,\\
&\Delta(x)=(-x)\otimes g+1\otimes (-x),\\
&\varepsilon(1_{H})=1,\ \varepsilon(g)=1,\ \varepsilon(x)=0,\\
&S(1_{H})=1_{H},\ S(g)=g,\ S(x)=-gx.
\end{align*}
Then $(H,\alpha)$ is a monoidal Hom-Hopf algebra.
Let $R=\frac{1}{2}(1\otimes 1+1\otimes g+g\otimes 1-g\otimes g)$. It is easy to check that $(H,\alpha,R)$ is quasitriangular.

\noindent{\bf Proposition 4.10.}
Let $(H,\alpha)$ be a quasitriangular monoidal Hom-Hopf algebra with the quasitriangular structure $R$. Then $R$ satisfies the QHYBEs
\begin{align}
&R^{12}(R^{13}R^{23})=(R^{13}R^{23})R^{12}\\
&(R^{12}R^{13})R^{23}=R^{23}(R^{13}R^{12}).
\end{align}

\noindent{\bf Proof}~~
The proof is essentially the same as in the Hopf algebra setting.$\hfill \Box$
\\

\noindent{\bf Proposition 4.11.}  The Drinfel'd double $(H^{op}\bowtie H^{\ast},\alpha\otimes(\alpha^{-1})^{\ast})$ has the quasitriangular structure
$$R=\sum\limits_{i} 1\ltimes h^{*}_{i}\o S^{-1}(h_{i})\ltimes \varepsilon,$$
where $\{h_{i}\}$ and $\{h^{*}_{i}\}$ are a base of $H$ and its dual base in $H^{\ast}$ respectively.

\noindent{\bf Proof}~~Firstly, obviously $(\varepsilon\otimes id)R=(id\otimes\varepsilon)R=1$.

Then for any $h\ltimes f\in H^{op}\bowtie H^{\ast}$, and $l,g\in H^{\ast}$, $a,b\in H$, on one hand, we have
$$\begin{aligned}
&\langle \Delta^{op}(h\ltimes f)R,g\otimes a\otimes l\otimes b\rangle\\
=&\sum\langle (h_{2}\ltimes f_{2}\otimes h_{1}\ltimes f_{1})R,g\otimes a\otimes l\otimes b\rangle\\
=&\sum\langle \alpha(h_{2})\ltimes f_{2}h^{\ast}_{i}\otimes\alpha^{2}(S^{-1}(h_{i})_{21})h_{1}\ltimes f_{1}(S(S^{-1}(h_{i})_{1})(\alpha^{-3}(?)S^{-1}(h_{i})_{22})),\\
&g\otimes a\otimes l\otimes b\rangle\\
=&\sum\langle g,\alpha(h_{2})\rangle\langle f_{2}h^{\ast}_{i},a\rangle\langle l,\alpha^{2}(S^{-1}(h_{i})_{21})h_{1}\rangle\langle f_{1},S(S^{-1}(h_{i})_{1})(\alpha^{-3}(b)S^{-1}(h_{i})_{22})\rangle\\
=&\sum\langle g,\alpha(h_{2})\rangle\langle f_{2},a_{1}\rangle\langle h^{\ast}_{i},a_{2}\rangle\langle l,\alpha^{2}(S^{-1}(h_{i})_{21})h_{1}\rangle\langle f_{1},S(S^{-1}(h_{i})_{1})(\alpha^{-3}(b)S^{-1}(h_{i})_{22})\rangle\\
=&\sum\langle g,\alpha(h_{2})\rangle\langle l,\alpha^{2}(S^{-1}(a_{212}))h_{1}\rangle\langle f,[a_{22}(\alpha^{-3}(b)S^{-1}(a_{211}))]a_{1}\rangle\\
=&\sum\langle g,\alpha(h_{2})\rangle\langle l,\alpha^{2}(S^{-1}(a_{212}))h_{1}\rangle\langle f,[a_{22}\alpha^{-2}(b)][S^{-1}\alpha(a_{211})\alpha^{-1}(a_{1})]\rangle\\
=&\sum\langle g,\alpha(h_{2})\rangle\langle l,S^{-1}(a_{1})h_{1}\rangle\langle f,a_{2}\alpha^{-1}(b)\rangle,
\end{aligned}$$
and on the other hand
$$\begin{aligned}
&\langle R\Delta(h\ltimes f),g\otimes a\otimes l\otimes b\rangle\\
=&\sum\langle (1\ltimes h^{\ast}_{i})(h_{1}\ltimes f_{1})\otimes(S^{-1}(h_{i})\ltimes\varepsilon_{H})(h_{2}\ltimes f_{2}),g\otimes a\otimes l\otimes b\rangle\\
=&\sum\langle \alpha^{3}(h_{121})\ltimes h^{\ast}_{i}(S(h_{11})(\alpha^{-2}(?)h_{122}))f_{1}\otimes h_{2}S^{-1}(h_{i})\ltimes f_{2}\circ\alpha^{-1},g\otimes a\otimes l\otimes b\rangle\\
=&\sum\langle g,\alpha^{3}(h_{121})\rangle\langle h^{\ast}_{i},S(h_{11})(\alpha^{-2}(a_{1})h_{122})\rangle\langle f_{1},a_{2}\rangle\langle l, h_{2}S^{-1}(h_{i})\rangle\langle f_{2},\alpha^{-1}(b)\rangle\\
=&\sum\langle g,\alpha^{3}(h_{121})\rangle\langle f,a_{2}\alpha^{-1}(b)\rangle\langle l, h_{2}[(S^{-1}(h_{122})S^{-1}(\alpha^{-2}(a_{1})))h_{11}]\rangle\\
=&\sum\langle g,\alpha^{2}(h_{21})\rangle\langle f,a_{2}\alpha^{-1}(b)\rangle\langle l, (\alpha(h_{222})S^{-1}\alpha(h_{221}))(S^{-1}\alpha^{-1}(a_{1})\alpha^{-1}(h_{1}))\rangle\\
=&\sum\langle g,\alpha(h_{2})\rangle\langle f,a_{2}\alpha^{-1}(b)\rangle\langle l, (S^{-1}(a_{1})h_{1}\rangle,
\end{aligned}$$
Hence $R\Delta(h\ltimes f)=\Delta^{op}(h\ltimes f)R.$

Similarly we have $(\Delta\otimes\alpha^{-1})R=R^{13}R^{23},$
$(\alpha^{-1}\otimes\Delta)R=R^{13}R^{12}.$ This completes the proof.$\hfill \Box$
\\

\noindent{\bf Example 4.12.}
In the Example 4.6, the Drinfel'd double of $(kG,\phi)$ is given where $G$ is a finite group and $\phi$ is a group automorphism of $G$. Then by the above  proposition, the quasitriangular structure of $D(kG)$ is
$$R=\sum_{g\in G}1_{G}\ltimes e_{g}\o g^{-1}\ltimes 1_{k^{G}}.$$

\section*{ACKNOWLEDGEMENTS}

This work was supported by the NSF of China (No. 11901240, 11801304) and the NSF of Shandong Province (No. ZR2018PA006).

\end{document}